\documentclass[12pt]{article}
\usepackage{verbatim}
\usepackage{latexsym}
\usepackage{amsmath,amsthm,amsfonts}
\usepackage{eucal}
\usepackage{cases}

\usepackage{mathtext}
\usepackage[cp1251]{inputenc}
\usepackage[T2A]{fontenc}
\usepackage[dvips]{graphicx}
\usepackage{amsmath}
\usepackage{amssymb}
\usepackage{amsxtra}
\usepackage{latexsym}
\usepackage{ifthen}
\usepackage{color}

\newtheorem{theorem}{Theorem}
\newtheorem{definition}{Definition}[section]

\numberwithin{equation}{section}

\begin{document}

\title{Strong solutions of the thin film equation in spherical geometry.}

\author{Roman M. Taranets \\
 Institute of Applied Mathematics and Mechanics of the NASU,\\
 Sloviansk, Ukraine, taranets\_r@yahoo.com}

\maketitle
\begin{abstract}
We study existence and long-time behaviour of strong solutions for
the thin film equation using a priori estimates in a weighted Sobolev space. 
This equation can be classified as a doubly degenerate fourth-order parabolic 
and it models coating flow on the outer surface of a sphere. It is shown that
the strong solution asymptotically decays to the flat profile.
\end{abstract}

\section{Introduction}

In this paper, we study the following doubly degenerate fourth-order parabolic equation
\begin{equation}\label{eq:main eqn}
u_{t}+\left( (1-x^{2} ) |u|^{n}( (1-x^{2} )u_{x} )_{xx}\right)_{x}=0 \text{ in } Q_T,
\end{equation}
where $Q_T = \Omega \times(0,T )$, $n > 0$, $T > 0$, and $\Omega=(-1,1)$.
This equation describes the dynamics of a thin viscous liquid film on the
outer surface of a solid sphere. More general dynamics of the liquid film for the cases
when the draining of the film due to gravity were balanced by centrifugal
forces arising from the rotation of the sphere about a vertical axis and by capillary forces due
to surface tension was considered in \cite{kangcoatingsphere}. In addition, Marangoni effects due to temperature gradients 
were taken into account in \cite{KNC-17}.
The spherical model without the surface tension and Marangoni effects was studied in \cite{TH-10,Wil-94}.

In \cite{kangcoatingsphere}, the authors derived the following equation for no-slip regime in dimensionless form
$$
h_{t}+ \tfrac{1}{\sin\theta} (h^{3} \sin  \theta  \, J  )_{\theta}=0,
$$
$$
J :=  a\sin\theta+b\sin\theta\cos\theta + c [ 2h+\tfrac{1}{\sin\theta} (\sin\theta \, h_{\theta} )_{\theta} ]_{\theta} ,
$$
where $h(\theta,t)$ represent the thickness of the thin film, $\theta \in (0,\pi)$ is the polar angle in spherical
coordinates, with $t$ denoting time; the dimensionless parameters $a$, $b$ and $c$ describe
the effects of gravity, rotation and surface tension, respectively.
After the change of variable $x=-\cos\theta$, this equation can be written in the form:
\begin{equation}\label{eq:physical eqn}
u_{t}+ [u^{3} (1-x^{2} ) (a - b x + c (2h + ((1-x^{2})u_{x} )_{x} )_{x} ) ]_{x}=0,
\end{equation}
where $x \in (-1,1)$. As a result, equation (\ref{eq:main eqn}) for $n =3$ is a particular case of (\ref{eq:physical eqn})
for no-slip regime. On the other hand, (\ref{eq:main eqn}) for $n < 3$ generalises (\ref{eq:physical eqn}) with
$a = b =0$ for different slip regimes, for example, like weak or partial wetting.

In contrast to the classical thin film equation:
\begin{equation}\label{eq:Bernis eqn}
u_{t}+\left(|u|^{n}u_{xxx}\right)_{x}=0,
\end{equation}
which describes the behavior of a thin viscous film on a flat surface under
the effect of surface tension, the equation (\ref{eq:main eqn}) is not yet well analysed.
To the best of our knowledge, there is only one analytical result \cite{KSZ-16} where the authors
proved existence of weak solutions in a weighted Sobolev space.  In 1990, Bernis and Friedman \cite{bernis1990higher}
constructed non-negative weak solutions of the equation (\ref{eq:Bernis eqn}) when
$n \geqslant 1$, and it was also shown that for $n\geqslant 4$, with a positive
initial condition, there exists a unique positive classical solution.
In 1994, Bertozzi et al. \cite{bertozzi1994singularities} generalised this positivity
property  for the case $n \geqslant \frac{7}{2}$. In 1995, Beretta et al. \cite{beretta1995nonnegative}
proved the existence of non-negative weak solutions for the equation (\ref{eq:Bernis eqn}) if
$n > 0$, and the existence of strong ones for $0< n < 3$. Also, they could show that this positivity-preserving property holds for almost every
time $t$ in the case $n \geqslant 2$.  A similar result on a cylindrical
surface was obtained in \cite{Ch10}. Regarding the long-time behaviour, Carrillo
and Toscani \cite{carrillo2002long} proved the convergence to a self-similar
solution for equation (\ref{eq:Bernis eqn}) with $n=1$ and Carlen
and Ulusoy \cite{carlen2007asymptotic} gave an upper bound on the
distance from the self-similar solution. A similar result on a cylindrical
surface was obtained in \cite{burchard2012convergence}.

In the present article, we obtain the existence of weak solutions in
a wider weighted classes of functions than it was done in \cite{KSZ-16}. Moreover, we show the
existence of non-negative strong solutions and we also prove that this solution decays 
asymptotically to the flat profile. Note that (\ref{eq:main eqn}) loses
its parabolicity not only at $u=0$ (as in (\ref{eq:Bernis eqn})) but also at $x=\pm1$.
For this reason, it is natural to seek solution in a Soblev space with weight $1-x^2$.
For example, it is the well-known that the non-negative steady state of equation (\ref{eq:Bernis eqn})
for $x \in (-1,1)$ has the form
$$
u_s(x) = c_1(1-x^2) + c_2, \text{ where } c_i \geqslant 0.
$$
On the other hand, the equation (\ref{eq:main eqn}) has the following non-negative steady state
$$
u_s(x) = (c_1 +c_2) \ln (1+x) + (c_1 - c_2) \ln (1- x) + c_3,
$$
where $0 \leqslant |c_2| \leqslant - c_1$,
$c_3 \geqslant - (c_1 +c_2) \ln (1+\frac{c_2}{c_1}) + (c_1 - c_2) \ln (1- \frac{c_2}{c_1})$,
hence $u_s(x) \to +\infty$ as $x \to \pm 1$.

\section{Existence of Strong Solutions}

We study the following thin film equation
\begin{equation}\label{A-1}
u_{t}+\left((1-x^{2}) |u|^{n}\left((1-x^{2})u_{x}\right)_{xx}\right)_{x}=0 \text{ in } Q_T
\end{equation}
with the no-flux boundary conditions
\begin{equation}\label{A-2}
(1-x^{2})u_{x}=(1-x^{2})\left((1-x^{2})u_{x}\right)_{xx}=0 \text{ at } x=\pm1,\, t > 0,
\end{equation}
and the initial condition
\begin{equation}\label{A-3}
u(x,0)=u_{0}(x).
\end{equation}
Here $n > 0$, $Q_T = \Omega \times (0,T)$, $\Omega :=(-1,1)$,  and $T >0$.
Integrating the equation (\ref{A-1}) by using boundary conditions (\ref{A-2}), we obtain the mass conservation property
\begin{equation}\label{mass-con}
\int \limits_{\Omega}{ u(x,t) dx} = \int \limits_{\Omega}{ u_0(x ) dx} =: M.
\end{equation}
Consider  initial data $u_0(x) \geqslant 0$ for all $x \in \bar{\Omega}$ satisfying
\begin{equation}\label{A-4}
\int \limits_{\Omega}{ \{(1-x^2)^{\beta} u^2_0(x) + (1-x^2) u^2_{0,x} (x)\} dx} < \infty,\ \beta \in (0,\tfrac{2}{n}] .
\end{equation}

\begin{definition}\label{def}[weak solution]
Let $n > 0$. A  function $u$ is a weak solution of
the problem (\ref{A-1})--(\ref{A-3}) with initial data $u_0$ satisfying (\ref{A-4})
if   $u(x,t)$ 
has the following properties
$$
(1-x^2)^{\frac{\beta}{2}} u \in C_{x,t}^{\frac{\alpha}{2}, \frac{\alpha}{8}}(\bar{Q}_T),
\ 0 < \alpha < \beta \leqslant \tfrac{2}{n} ,
$$
$$
u_t \in L^2(0,T; (H^1(\Omega))^*), \ (1-x^2)^{\frac{1}{2}}  u_x  \in L^{\infty}(0,T; L^2(\Omega)),
$$
$$
(1-x^2)^{\frac{1}{2}} |u|^{\frac{n}{2}} ( (1-x^2)u_x )_{xx}   \in L^2(P),
$$
$u$ satisfies (\ref{A-1}) in the following sense:
$$
\int \limits_{0}^T {\langle u_t , \phi \rangle \,dt} -
\iint \limits_{P} { (1-x^2) |u|^{n} ( (1-x^2)u_x )_{xx} \phi_{x} \,dx dt}\\
 =0
$$
for all $\phi \in L^2(0,T; H^1(\Omega))$, where $P : = \bar{Q}_T \setminus
\{ \{u = 0\} \cup \{t =0\}\} $,
$$
(1-x^2)^{\frac{1}{2}} u_x(., t) \to (1-x^2)^{\frac{1}{2}} u_{0,x}(.)
\text{ strongly in } L^2(\Omega) \text{ as } t \to 0,
$$
and boundary conditions (\ref{A-2}) hold at all points of the lateral boundary, where
$\{u > 0\}$.
\end{definition}

Let us denote by
\begin{equation}\label{G_0-def}
0 \leqslant G_{0}(z): =   \begin{cases} \tfrac{ z^{2 -n} - A^{2 -n} }{(n-1)(n
- 2)} - \tfrac{A^{1-n}}{1 -n}(z - A) \text{
if } n \neq   1, 2  , \\
z \ln z - (z - A)(\ln A + 1)  \text{ if }
n = 1 , \\
\ln (\tfrac{A}{z}) + \tfrac{z}{A} - 1 \text{ if } n = 2,
\end{cases}
\end{equation}
where $A = 0$ if $n \in (1,2)$ and $A > 0$ if else. Next, we establish existence of a more regular
solution $u$ of (\ref{A-1}) than a weak solution in  the sense of Definition~\ref{def}. Besides,
we show that this strong solution $u$ with some weight exponentially decays to zero.

\begin{theorem}[strong solution]\label{Th1}
Assume that $n \geqslant 1$ and initial data $u_0$ satisfies $\int \limits_{\Omega}{G_0(u_0) \,dx} < +\infty $ then the problem (\ref{A-1})--(\ref{A-3}) has a non-negative weak solution, $u$, in the sense of Definition~\ref{def},
such that
$$
(1-x^2)u_x  \in L^2(0,T;H^1(\Omega)), \ (1-x^2)^{\frac{\gamma}{2}}u_x  \in L^2(Q_T), \ \gamma \in (0,1).
$$
$$
(1-x^2)^{\frac{\mu}{2}}u   \in L^2(Q_T),  \  \mu \in (-1, \beta).
$$
Moreover, there exist positive constants $A$, $B$ depending on initial data such
that
$$
\tfrac{1}{2} \int \limits_{\Omega} {  (1-x^{2} ) u^2_{ x}(x,t)   \,  dx} \leqslant A \, e^{- B\,t} \
\forall \, t \geqslant 0,
$$
hence
$$
(1-x^2)^{\frac{\beta}{2}} |u  - \tfrac{M}{|\Omega|}|  \to 0 \text{ as } t \to +\infty.
$$
\end{theorem}

\section{Proof of Theorem~\ref{Th1}}

\subsection{Approximating problems}

Equation (\ref{A-1}) is doubly degenerate when $u =0$ and $x = \pm 1$. For this reason, for
any $\epsilon > 0$ and $\delta > 0$ we consider two-parametric regularised equations
\begin{equation}\label{eq:regularized eqn}
 u_{\epsilon \delta,t} + \left[(1-x^{2}+\delta) (|u_{\epsilon\delta}|^{n}+ \epsilon  )\left((1-x^{2}+\delta) u_{\epsilon \delta,x}\right)_{xx}\right]_{x}=0 \text{ in } Q_T
\end{equation}
with boundary conditions
\begin{equation}\label{reg-1}
 u_{\epsilon \delta,x} =\left((1-x^{2}+\delta) u_{\epsilon \delta,x} \right)_{xx}=0  \text{ at } x= \pm1,
\end{equation}
and initial data
\begin{equation}\label{reg-2}
u_{\epsilon \delta}(x,0)= u_{0,\epsilon \delta }(x)\in C^{4 + \gamma}( \bar{\Omega}), \ \gamma > 0,
\end{equation}
where
\begin{equation}\label{reg-3-0}
u_{0,\epsilon\delta }(x) \geqslant u_{0\delta }(x) + \epsilon^{\theta}, \  \ \theta \in (0, \tfrac{1}{2(n-1)}),
\end{equation}
\begin{equation}\label{reg-3-1}
  u_{0 ,\epsilon \delta} \to   u_{0 \delta} \text{ strongly in } H^1(\Omega)
\text{ as } \epsilon \to 0,
\end{equation}
\begin{equation}\label{reg-3}
(1-x^{2}+\delta)^{\frac{1}{2}} u_{0x,\delta} \to (1-x^{2}) u_{0,x } \text{ strongly in } L^2(\Omega)
\text{ as }  \delta \to 0.
\end{equation}
The parameters $\epsilon > 0$ and $\delta > 0$ in (\ref{eq:regularized eqn}) make  the problem
regular up to the boundary (i.e. uniformly parabolic).
The existence of a solution of (\ref{eq:regularized eqn}) in a small
time interval is guaranteed by the Schauder estimates in \cite{friedman1958interior}.
Now suppose that $u_{\epsilon \delta}$ is a solution of equation
(\ref{eq:regularized eqn}) and that it is continuously differentiable
with respect to the time variable and fourth order continuously differentiable
with respect to the spatial variable.

\subsection{Existence of weak solutions}

In order to get an \textit{a priori} estimation of $u_{\epsilon \delta}$, we multiply both sides
of equation (\ref{eq:regularized eqn}) by $- [(1-x^{2}+\delta) u_{\epsilon \delta,x}]_{x}$
and integrate over $\Omega$ by (\ref{reg-1}). This gives us
\begin{multline}\label{eq:estimate IBP-0}
\tfrac{1}{2} \tfrac{d}{dt} \int \limits_{\Omega}{ (1-x^{2}+\delta) u_{\epsilon \delta, x}^2 \, dx}  +     \\
\int \limits_{\Omega}{(1-x^{2}+\delta) (|u_{\epsilon \delta}|^{n}+\epsilon )
[(1-x^{2}+\delta)u_{\epsilon \delta,x} ]_{xx}^2\,dx} =0.
  \end{multline}
Integrating (\ref{eq:estimate IBP-0}) in time, we get
\begin{multline} \label{eq:estimate IBP}
\tfrac{1}{2}  \int \limits_{\Omega}{ (1-x^{2}+\delta) u_{\epsilon \delta, x}^2 \, dx} +  \\
 \iint \limits_{Q_T}{(1-x^{2}+\delta) (|u_{\epsilon \delta}|^{n}+\epsilon )
[(1-x^{2}+\delta)u_{\epsilon \delta,x} ]_{xx}^2\,dx dt} = \\
\tfrac{1}{2}
\int \limits_{\Omega}{ (1-x^{2}+\delta) u_{0x,\epsilon \delta}^2 \, dx}  .
\end{multline}
By (\ref{reg-3}) we have
\begin{equation}\label{eq:estimate ux bound}
\int \limits_{\Omega}{ (1-x^{2}+\delta) u_{\epsilon \delta, x}^2 \, dx} \leqslant C_{0},
\end{equation}
where $C_{0}>0$ is independent of $\epsilon$ and $\delta$. From (\ref{eq:estimate ux bound})
and (\ref{eq:estimate IBP}) it follows that
\begin{equation}\label{ap-1}
\{u_{\epsilon \delta}\}_{\epsilon > 0} \text{ is uniformly bounded in } L^{\infty}(0,T;H^1(\Omega)).
\end{equation}
\begin{equation}\label{ap-2}
\{(1-x^{2}+\delta)^{\frac{1}{2}} (|u_{\epsilon \delta}|^{n}+\epsilon )^{\frac{1}{2}}
[(1-x^{2}+\delta)u_{\epsilon \delta,x} ]_{xx} \}_{\epsilon,\,\delta > 0} \text{ is u.\,b. in } L^{2}(Q_T).
\end{equation}
By (\ref{ap-1}) and (\ref{ap-2}), using the same method as \cite{bernis1990higher}, we can
prove that solutions $u_{\epsilon \delta}$ have uniformly (in $\epsilon$) bounded  $C_{x,t}^{1/2,1/8}$-norms.
By the Arzel\`{a}-Ascoli theorem, this equicontinuous property, together with the uniformly boundedness
shows that  every sequence $\{u_{\epsilon \delta}\}_{\epsilon > 0}$ has
a subsequence   such that
\begin{equation}\label{eq:uniformly convergence}
u_{\epsilon \delta} \to u_{\delta} \text{ uniformly in } Q_T \text{ as }\epsilon \to 0.
\end{equation}
As a result, we obtain a solution $u_{\delta}$ of the problem (\ref{eq:regularized eqn})--(\ref{reg-2}) with
$\epsilon = 0$ in the sense of \cite[Theorem~3.1, pp. 185--186]{bernis1990higher}.

Next, we show that the family of solutions $\{u_{ \delta}\}_{\delta > 0}$ is uniformly bounded in
some weighted space. Using the mass conservation property
$$
\int \limits_{\Omega}{ u_{\delta}(x,t) dx} =   M > 0,
$$
we arrive at
\begin{equation}\label{est-00}
|u_{\delta} - \tfrac{M}{|\Omega|}| = \Bigl | \int \limits_{x_0}^{x} { u_x \,dx } \Bigr |
\leqslant \Bigl( \int \limits_{\Omega}  { (1- x^2) u^2_x \,dx } \Bigr)^{\frac{1}{2}}
\Bigl( \int \limits_{x_0}^{x}  {\tfrac{dx}{1- x^2}  } \Bigr)^{\frac{1}{2}}.
\end{equation}
Multiplying (\ref{est-00}) by $(1-x^2)^{\frac{\beta}{2}}$ for any $\beta > 0$, by (\ref{eq:estimate ux bound})
we deduce that
\begin{equation}\label{est-001}
(1-x^2)^{\frac{\beta}{2}} |u_{\delta} - \tfrac{M}{|\Omega|}| \leqslant (\tfrac{C_0}{2})^{\frac{1}{2}}
\Bigl( (1-x^2)^{ \beta } \ln (\tfrac{(1+x)(1-x_0)}{(1-x)(1+x_0)}) \Bigr)^{\frac{1}{2}} \leqslant C_1
\end{equation}
for all $x \in \bar{\Omega}$, where $C_1 > 0$ is independent of $\delta > 0$. From (\ref{est-001}) we find that
\begin{equation}\label{ap-4}
\{(1-x^{2} )^{\frac{\beta}{2}} u_{ \delta } \}_{ \delta > 0} \text{ is u.\,b. in } Q_T \text{ for any } \beta > 0.
\end{equation}
In particular, by (\ref{eq:estimate ux bound}) we get
\begin{equation}\label{eq:equivcont x}
 (1-x^{2})^{\frac{\beta}{2}} | u_{ \delta}(x_{1},t) - u_{ \delta}(x_{2},t) |
 \leqslant C_2 |x_{1}-x_{2}|^{ \frac{\alpha}{2}} \ \ \forall \, x_{1},x_{2} \in  \Omega , \ \alpha \in (0,\beta).
\end{equation}
By (\ref{ap-2}), (\ref{ap-4}) and (\ref{eq:equivcont x}) with $\beta \in (0, \frac{2}{n}]$, using the same method as \cite[Lemma 2.1, p.183]{bernis1990higher}, we can
prove similarly that
\begin{equation}\label{eq:equivcont t}
(1-x^{2})^{\frac{\beta}{2}} | u_{ \delta}(x,t_{1}) - u_{ \delta}(x,t_{2}) | \leqslant C_{3}|t_{1}-t_{2}|^{ \frac{\alpha}{8}} \ \ \forall\, t_{1},t_{2}\in(0,T).
\end{equation}
The inequalities (\ref{eq:equivcont x}) and (\ref{eq:equivcont t})
show the uniform (in $\delta$) boundedness of a sequence $\{ (1-x^{2})^{\frac{\beta}{2}} u_{ \delta} \}_{\delta > 0}$
in the $C_{x,t}^{\frac{\alpha}{2}, \frac{\alpha}{8}} $-norm.
By the Arzel\`{a}-Ascoli theorem, this a priori bound together with (\ref{ap-4})
 shows that as $\delta \to 0$, every sequence $\{(1-x^{2})^{\frac{\beta}{2}} u_{ \delta}\}_{\delta > 0}$ has
a subsequence $\{(1-x^{2})^{\frac{\beta}{2}} u_{ \delta_k}\}_{\delta_k > 0}$ such that
\begin{equation}\label{eq:uniformly convergence}
(1-x^{2})^{\frac{\beta}{2}} u_{ \delta_k} \to  (1-x^{2})^{\frac{\beta}{2}} u
\text{ uniformly in } \bar{Q}_T  \text{ as } \delta_k \to 0.
\end{equation}
Following the idea of proof \cite[Theorem~3.1]{bernis1990higher}, we  obtain a solution $u $ of the problem
(\ref{eq:regularized eqn})--(\ref{reg-2}) in the sense of Definition~\ref{def}.

\subsection{Existence of strong solutions}

Let us denote by $G_{\epsilon}(z)$ the following function
$$
G_{\epsilon}(z) \geqslant 0 \ \forall\, z \in \mathbb{R}, \ G''_{\epsilon}(z)= \tfrac{1}{|s|^{n}+\epsilon}.
$$
Now we multiply equation (\ref{eq:regularized eqn})
by $G'_{\epsilon}(u_{\epsilon \delta})$ and integrate over $\Omega$
to get
\begin{equation}\label{eq:positivity proof integrate}
\tfrac{d}{dt} \int \limits_{\Omega}{ G_{\epsilon} (u_{\epsilon \delta}(x,t) ) \, dx} +
\int \limits_{\Omega} { [(1-x^{2}+\delta) u_{\epsilon \delta,x} ]_{ x}^2 \,  dx} =0.
\end{equation}
After  integration in time, equation (\ref{eq:positivity proof integrate})
becomes
\begin{equation}\label{entr-1}
 \int \limits_{\Omega}{ G_{\epsilon} (u_{\epsilon \delta}(x,T) ) \, dx} +
\iint \limits_{Q_T} { [(1-x^{2}+\delta) u_{\epsilon \delta,x} ]_{ x}^2 \,  dxdt} =
\int \limits_{\Omega}{ G_{\epsilon} (u_{0,\epsilon \delta }(x ) ) \, dx}.
\end{equation}
We compute
$$
G''_0(z) - G''_{\epsilon}(z) = \tfrac{\epsilon}{|z|^n(|z|^n + \epsilon)},
$$
and consequently
$$
G_0(z) - G_{\epsilon}(z) = \epsilon \int \limits_A^z { \int \limits_A^v { \tfrac{ds dv}{|s|^n(|s|^n + \epsilon)} } },
$$
where $A$ is some positive constant. As $u_{0, \epsilon\delta}(x )$ is bounded then by (\ref{reg-3-0}) it follows that
$$
|G_0(u_{0, \epsilon \delta}(x )) - G_{\epsilon}(u_{0, \epsilon\delta }(x ))| \leqslant C\,\epsilon^{1 - 2 \theta(n-1)} \to 0
\text{ as } \epsilon \to 0,
$$
and therefore, due to (\ref{reg-3-1}), we have
\begin{equation}\label{entr-2}
\int \limits_{\Omega}{ G_{\epsilon} (u_{0,\epsilon  }(x ) ) \, dx} \to \int \limits_{\Omega}{ G_{0} (u_{0 \delta}(x ) ) \, dx}
\text{ as } \epsilon \to 0.
\end{equation}
As a result, by (\ref{entr-1}), (\ref{entr-2}) we deduce that
\begin{equation}\label{eq:positivity proof upper bound}
\int \limits_{\Omega} {G_{\epsilon} (u_{\epsilon \delta}(x,T) ) \, dx}  \leqslant C_4,
\end{equation}
\begin{equation}\label{ap-3}
\{(1-x^{2}+\delta) u_{\epsilon \delta, x} \}_{\epsilon,\, \delta > 0} \text{ is u.\,b. in } L^{2}(0,T;H^1(\Omega)),
\end{equation}
where $C_1 > 0$ is independent of $\epsilon$ and $\delta$.
Similar to \cite[Theorem~4.1, p. 190]{bernis1990higher}, using
(\ref{ap-1}) and (\ref{eq:positivity proof upper bound}),
we can show that the limit solution $ u_{ \delta} $ is non-negative if $n \in [1,4)$ and
positive if $n \geqslant 4$. Next, letting $\delta \to 0$, we get a non-negative strong
solution.

\subsection{Asymptotic behaviour}

Let us denote  by
$$
\mathcal{E}_{\delta}(u(t)) : = \tfrac{1}{2} \int \limits_{\Omega} {  (1-x^{2} + \delta) u^2_{ x}   \,  dx}, \  \
\tfrac{\delta \mathcal{E}_{\delta}(u )}{\delta u} := - ( (1-x^{2} + \delta) u_{ x})_x.
$$
By using the notations, we rewrite (\ref{eq:estimate IBP-0}) and (\ref{eq:positivity proof integrate})
with $\epsilon =0$ in the form
\begin{equation} \label{ineq-1}
\tfrac{d}{dt} \mathcal{E}_{\delta} (u_\delta(t)) +   \int \limits_{\Omega}{(1-x^{2} + \delta ) u_{\delta}^{n}
[\tfrac{\delta \mathcal{E}_{\delta}(u_{\delta})}{\delta u}]_x^2 \,dx} =0 ,
\end{equation}
\begin{equation}\label{ineq-2}
\tfrac{d}{dt} \int \limits_{\Omega}{ G_0 (u_{\delta} ) \, dx} +
\int \limits_{\Omega} { [\tfrac{\delta \mathcal{E}_{\delta}(u_{\delta} )}{\delta u}]^2 \,  dx} =0.
\end{equation}
Next, we will use the following Hardy's inequality
\begin{equation}\label{hardy}
\int \limits_{-1}^1{(1-x^{2} )^{  -2 + \gamma} v^2(x) \,dx} \leqslant C
\int \limits_{-1}^1{ (1-x^2)^{ \gamma } v_x^2(x) \,dx}
\end{equation}
for any $\gamma > 0$ and for all $v \in H^1(-1,1)$ such that $v(\pm 1) = 0$. Really, using
integration by parts and Cauchy inequality, we have
\begin{multline*}
\int \limits_{-1}^1 {(1-x^{2} )^{  -2 + \gamma} v^2(x) \,dx}  = v^2(x)g(x) \biggl |_{-1}^1 -
2 \int \limits_{-1}^1 {  v(x) v_x(x) g(x) \,dx} \leqslant    \\
2
\Bigl( \int \limits_{-1}^1 { (1-x^{2} )^{ -2 +\gamma} v^2(x) \,dx} \Bigr)^{\frac{1}{2}}
\Bigl( \int \limits_{-1}^1 { (1-x^{2} )^{2 -\gamma }g^2(x) v_x^2(x) \,dx} \Bigr)^{\frac{1}{2}} ,
\end{multline*}
where
\begin{multline*}
|(1-x^2)^{1 - \frac{\gamma}{2}} g(x)| =
\Bigl | (1-x^2)^{1 - \frac{\gamma}{2}} \int \limits^x { (1-x^2)^{  - 2 + \gamma} dx } \Bigr |  \leqslant \\
 C\, (1-x^2)^{\frac{\gamma}{2}} \ \forall \, x \in [-1,1], \ \gamma \geqslant 0,
\end{multline*}
and as $v(x) \in C^{\frac{1}{2}}[-1,1]$ and $g(x) \sim (1-x^2)^{-1 + \gamma}$ at $x = \pm 1$
then $v^2(x)g(x)  = 0$ at $x = \pm 1$. From here we find that
\begin{multline*}
 \int \limits_{-1}^1 {(1-x^{2} )^{  -2 +\gamma } v^2(x) \,dx} \leqslant   \\
C
\Bigl( \int \limits_{-1}^1 { (1-x^{2} )^{  -2 +\gamma } v^2(x) \,dx} \Bigr)^{\frac{1}{2}}
\Bigl( \int \limits_{-1}^1 { (1-x^{2} )^{ \gamma }  v_x^2(x) \,dx} \Bigr)^{\frac{1}{2}} ,
\end{multline*}
whence it follows (\ref{hardy}).

Applying  (\ref{hardy}) to $v = (1-x^{2} ) u_{ x}$ with $\gamma = 1$, we obtain that
$$
\int \limits_{\Omega} {(1-x^{2} ) u^2_{ x}  \,dx} \leqslant C_5
\int \limits_{\Omega}  { (1-x^{2} ) [(1-x^{2} )  u_{ x}]_x^2 \,dx} \leqslant
C_5 \int \limits_{\Omega}  { [(1-x^{2} )  u_{ x}]_x^2 \,dx},
$$
i.\,e.
\begin{equation}\label{hardy-1}
2 \mathcal{E}_0(u(t)) \leqslant C_5 \int \limits_{\Omega} { [\tfrac{\delta \mathcal{E}_0(u )}{\delta u}]^2 \,  dx}.
\end{equation}
Summing (\ref{ineq-1}) and (\ref{ineq-2}), after integrating in time, taking $\delta \to 0$,
 and using (\ref{hardy-1}),
we arrive at
\begin{equation} \label{ineq-3}
 \mathcal{E}_0(u(t)) + B \int \limits_{0}^t { \mathcal{E}_0(u(s)) \,ds}
 \leqslant A: = \mathcal{E}(u_0) +  \int \limits_{\Omega}{ G_0 (u_0 ) \, dx} ,
\end{equation}
where $B:= \frac{2}{C_5}$. From (\ref{ineq-3}) by comparing to the solution $y(t)$ of the problem for ODE
$$
y'(t) + By(t) = 0, \ \ y(0) = A,
$$
we get
\begin{equation}\label{asym-1}
0 \leqslant \mathcal{E}_0(u(t)) \leqslant A \, e^{-B\, t} \to 0 \text{ as } t \to +\infty.
\end{equation}
By (\ref{est-001}) and (\ref{asym-1}) we deduce that
$$
(1-x^2)^{\frac{\beta}{2}} |u  - \tfrac{M}{|\Omega|}| \leqslant \tilde{A} \, e^{ -\tilde{B}\,t} \to 0
\text{ as } t \to +\infty.
$$
This proves Theorem~\ref{Th1} completely. $\square$


\begin{thebibliography}{99}

\bibitem{beretta1995nonnegative}
E. Beretta, M. Bertsch, and R. Dal Passo.
\newblock Nonnegative solutions of a fourth-order nonlinear degenerate parabolic equation.
\newblock {\em Archive for rational mechanics and analysis}, 129(2): 175--200, 1995.

\bibitem{bernis1990higher}
F. Bernis, A. Friedman.
\newblock Higher order nonlinear degenerate parabolic equations.
\newblock {\em J. Differential Equations}, 83(1): 179--206, 1990.

\bibitem{bertozzi1994singularities}
Andrea L. Bertozzi  et al.
\newblock {\em Singularities and similarities in interface flows.}
\newblock Trends and perspectives in applied mathematics. Springer New York, 155--208, 1994.


\bibitem{burchard2012convergence}
\newblock Almut  Burchard,  Marina Chugunova, and Benjamin K. Stephens.
\newblock Convergence to equilibrium for a thin-film equation on a cylindrical surface.
\newblock {\em  Communications in Partial Differential Equations}, 37(4): 585--609, 2012.


\bibitem{carlen2007asymptotic}
Eric A. Carlen, and S{\"u}leyman  Ulusoy.
 \newblock  Asymptotic equipartition and long time behavior of solutions of a thin-film equation.
 \newblock {\em Journal of Differential Equations}, 241(2): 279--292, 2007.

\bibitem{carrillo2002long}
Jos{\'e} A. Carrillo, and Giuseppe Toscani.
\newblock  Long-Time Asymptotics for Strong Solutions of the Thin Film Equation.
\newblock {\em Communications in mathematical physics}, 225(3): 551--571, 2002.


\bibitem{Ch10}
Marina Chugunova, Mary C. Pugh, and Roman M. Taranets.
\newblock Nonnegative solutions for a long-wave unstable
   thin film equation with convection.
\newblock {\em SIAM Journal on Mathematical Analysis}, 42(4): 1826--1853, 2010.

\bibitem{friedman1958interior}
Avner Friedman.
\newblock  Interior estimates for parabolic systems of partial differential equations.
\newblock {\em J. Math. Mech.}, 7(3): 393--417, 1958.


\bibitem{kangcoatingsphere}
D. Kang,  A. Nadim, and M. Chugunova.
\newblock  Dynamics and equilibria of thin viscous coating films on a rotating sphere.
\newblock {\em Journal of Fluid Mechanics}, 791: 495--518, 2016.

\bibitem{KNC-17}
D. Kang,  A. Nadim, and M. Chugunova.
\newblock  Marangoni effects on a thin liquid film coating a sphere with axial or radial thermal
gradients.
\newblock {\em Physics of Fluids}, 29: 072106-1--072106-15, 2017.

\bibitem{KSZ-16}
D. Kang,  Tharathep Sangsawang and Jialun Zhang.
\newblock  Weak solution of a doubly degenerate parabolic equation.
\newblock {\em  arXiv:1610.06303v2}, 2017.


\bibitem{TH-10}
D. Takagi, and Herbert E. Huppert.
\newblock  Flow and instability of thin films on a cylinder and sphere.
\newblock {\em Journal of Fluid Mechanics}, 647: 221--238, 2010.

\bibitem{Wil-94}
S.K. Wilson.
\newblock  The onset of steady Marangoni convection in a spherical geometry.
\newblock {\em Journal of Engineering Mathematics}, 28: 427--445, 1994.



\end{thebibliography}
\end{document}